%% file: Text-V.tex
\documentclass[11pt,twoside]{article}

\input{packages.tex}

\title{
  {\huge Complex Analysis of Real Functions \\[1.5ex]}
  V: The Dirichlet Problem on the Plane }

\author{
  \Large Jorge L. deLyra\footnote{Email: delyra@latt.if.usp.br} \\
  Department of Mathematical Physics \\
  Physics Institute \\
  University of São Paulo }

\date{March 19, 2018}


\input{math-defs.tex}

\begin{document}\maketitle

\vspace{-3.5ex}
\begin{abstract}
  \noindent
  In the context of the correspondence between real functions on the unit
  circle and inner analytic functions within the open unit disk, that was
  presented in previous papers, we show that the constructions used to
  establish that correspondence lead to very general proofs of existence
  of solutions of the Dirichlet problem on the plane. At first, this
  establishes the existence of solutions for almost arbitrary integrable
  real functions on the unit circle, including functions which are
  discontinuous and unbounded. The proof of existence is then generalized
  to a large class of non-integrable real functions on the unit circle.
  Further, the proof of existence is generalized to real functions on a
  large class of other boundaries on the plane, by means of conformal
  transformations.
\end{abstract}

\section{Introduction}\label{Sec01}

In previous papers~\cite{CAoRFI,CAoRFII,CAoRFIII,CAoRFIV} we have shown
that there is a correspondence between, on the one hand, integrable real
functions, singular Schwartz distributions and non-integrable real
functions which are locally integrable almost everywhere, and on the other
hand, inner analytic functions within the open unit disk of the complex
plane. This correspondence is based on the complex-analytic structure
within the unit disk of the complex plane, which we introduced
in~\cite{CAoRFI}. In order to establish this correspondence for integrable
and non-integrable real functions, we presented in~\cite{CAoRFI}
and~\cite{CAoRFIV} certain constructions which, given just such a real
function, produce from it a unique corresponding inner analytic function.

In this paper we will show that these constructions have, as collateral
consequences, the establishment of very general constructive proofs of the
existence of the solution of the Dirichlet boundary value problem for the
Laplace equation on regions of the plane. We will first establish the
proof for integrable real functions on the unit circle, then generalize it
to non-integrable real functions which are locally integrable almost
everywhere on that circle. Furthermore, with the use of conformal
transformations it is possible to generalize the proof to integrable and
non-integrable real function on other boundaries on the plane. We will
first establish this generalization for a large class of differentiable
curves, and then, with a single weak additional limitation on the real
functions, for a large class of curves that can be non-differentiable at a
finite set of points, such as polygons.

For ease of reference, we include here a one-page synopsis of the
complex-analytic structure introduced in~\cite{CAoRFI}. It consists of
certain elements within complex analysis~\cite{CVchurchill}, as well as of
their main properties.

\paragraph{Synopsis:} The Complex-Analytic Structure\\

\noindent
An {\em inner analytic function} $w(z)$ is simply a complex function which
is analytic within the open unit disk. An inner analytic function that has
the additional property that $w(0)=0$ is a {\em proper inner analytic
  function}. The {\em angular derivative} of an inner analytic function is
defined by

\noindent
\begin{equation}
  w^{\ldot}(z)
  =
  \ii
  z\,
  \frac{dw(z)}{dz}.
\end{equation}

\noindent
By construction we have that $w^{\ldot}(0)=0$, for all $w(z)$. The {\em
  angular primitive} of an inner analytic function is defined by

\begin{equation}
  w^{-1\ldot}(z)
  =
  -\ii
  \int_{0}^{z}dz'\,
  \frac{w(z')-w(0)}{z'}.
\end{equation}

\noindent
By construction we have that $w^{-1\ldot}(0)=0$, for all $w(z)$. In terms
of a system of polar coordinates $(\rho,\theta)$ on the complex plane,
these two analytic operations are equivalent to differentiation and
integration with respect to $\theta$, taken at constant $\rho$. These two
operations stay within the space of inner analytic functions, they also
stay within the space of proper inner analytic functions, and they are the
inverses of one another. Using these operations, and starting from any
proper inner analytic function $w^{0\ldot}(z)$, one constructs an infinite
{\em integral-differential chain} of proper inner analytic functions,

\begin{equation}
  \left\{
    \ldots,
    w^{-3\ldot}(z),
    w^{-2\ldot}(z),
    w^{-1\ldot}(z),
    w^{0\ldot}(z),
    w^{1\ldot}(z),
    w^{2\ldot}(z),
    w^{3\ldot}(z),
    \ldots\;
  \right\}.
\end{equation}

\noindent
Two different such integral-differential chains cannot ever intersect each
other. There is a {\em single} integral-differential chain of proper inner
analytic functions which is a constant chain, namely the null chain, in
which all members are the null function $w(z)\equiv 0$.

A general scheme for the classification of all possible singularities of
inner analytic functions is established. A singularity of an inner
analytic function $w(z)$ at a point $z_{1}$ on the unit circle is a {\em
  soft singularity} if the limit of $w(z)$ to that point exists and is
finite. Otherwise, it is a {\em hard singularity}. Angular integration
takes soft singularities to other soft singularities, and angular
differentiation takes hard singularities to other hard singularities.

Gradations of softness and hardness are then established. A hard
singularity that becomes a soft one by means of a single angular
integration is a {\em borderline hard} singularity, with degree of
hardness zero. The {\em degree of softness} of a soft singularity is the
number of angular differentiations that result in a borderline hard
singularity, and the {\em degree of hardness} of a hard singularity is the
number of angular integrations that result in a borderline hard
singularity. Singularities which are either soft or borderline hard are
integrable ones. Hard singularities which are not borderline hard are
non-integrable ones.

Given an integrable real function $f(\theta)$ on the unit circle, one can
construct from it a unique corresponding inner analytic function $w(z)$.
Real functions are obtained through the $\rho\to 1_{(-)}$ limit of the
real and imaginary parts of each such inner analytic function and, in
particular, the real function $f(\theta)$ is obtained from the real part
of $w(z)$ in this limit. The pair of real functions obtained from the real
and imaginary parts of one and the same inner analytic function are said
to be mutually Fourier-conjugate real functions.

Singularities of real functions can be classified in a way which is
analogous to the corresponding complex classification. Integrable real
functions are typically associated with inner analytic functions that have
singularities which are either soft or at most borderline hard. This ends
our synopsis.

\vspace{2.6ex}

\noindent
The treatment of the Dirichlet problem is usually developed under the
hypothesis that the boundary conditions are given by continuous real
functions at the boundary, leading to solutions which are continuous and
twice differentiable, with continuous derivatives, within the interior.
For the two-dimensional problems we will consider here, we will be able to
relax the conditions on the real functions at the boundary, accepting as
valid boundary conditions real functions which may not be continuous, and
not even bounded, at a finite set of boundary points. In order to allow
for this, the condition that the solution within the interior reproduces
the boundary condition everywhere at the boundary will have to be relaxed
to the reproduction only almost everywhere at the boundary. On the other
hand, it will also follow from the proofs offered that the solutions
within the interior are not only continuous and twice differentiable, but
in fact that they are always infinitely differentiable functions, on both
their arguments.

We begin our work in this paper in Section~\ref{Sec02}, by establishing
the existence theorem for boundary conditions given by integrable real
functions on the unit circle. This is followed, in Section~\ref{Sec03}, by
an extension of the existence theorem to non-integrable real functions on
the unit circle, which are, however, locally integrable almost everywhere
there. In Section~\ref{Sec04} we discuss the conformal transformations
that are required for the further versions of the existence theorem, that
are established in the subsequent sections. In Section~\ref{Sec05} we
establish the existence theorem for integrable real functions on almost
arbitrary differentiable simple closed curves on the plane. In
Section~\ref{Sec06} this existence theorem is extended to the case of
integrable real functions on non-differentiable simple closed curves on
the plane, curves which have, however, at most a finite set of points of
non-differentiability. In Section~\ref{Sec07} the existence theorem is
further extended, this time to non-integrable real functions, as qualified
above, on differentiable simple closed curves. Finally, in
Section~\ref{Sec08} we present the last and most general extension of the
existence theorem, to non-integrable real functions, as qualified above,
on non-differentiable simple closed curves, also as qualified above.

When we discuss real functions in this paper, some properties will be
globally assumed for these functions, just as was done in the previous
papers~\cite{CAoRFI,CAoRFII,CAoRFIII,CAoRFIV} leading to this one. These
are rather weak conditions to be imposed on these functions, that will be
in force throughout this paper. It is to be understood, without any need
for further comment, that these conditions are valid whenever real
functions appear in the arguments. These weak conditions certainly hold
for any real functions that are obtained as restrictions of corresponding
inner analytic functions to the unit circle, or to other simple closed
curves with finite total length.

The most basic global condition we will impose is that the real functions
must be measurable in the sense of Lebesgue, with the usual Lebesgue
measure~\cite{RARudin,RARoyden}. The second global condition we will
impose is that the real functions have no removable singularities. The
third and last global condition is that the number of hard singularities
of the real functions on their domains of definition be finite, and hence
that they be all isolated from one another. There will be no limitation on
the number of soft singularities.

The material contained in this paper is a development, reorganization and
extension of some of the material found, sometimes still in rather
rudimentary form, in the
papers~\cite{FTotCPI,FTotCPII,FTotCPIII,FTotCPIV,FTotCPV}.

\section{Integrable Real Functions on the Unit Circle}\label{Sec02}

In a previous paper~\cite{CAoRFI} we have shown that, given an integrable
real function on the unit circle, one can define from it a unique inner
analytic function whose real part reproduces that real function when
restricted to the unit circle. What follows is an outline of the
construction of this inner analytic function. Given the integrable real
function $f(\theta)$, we define from it, by means of the usual integrals,
the Fourier coefficients $\alpha_{0}$, $\alpha_{k}$ and $\beta_{k}$, for
$k\in\{1,2,3,\ldots,\infty\}$, and from those coefficients we define the
complex Taylor coefficients $c_{0}=\alpha_{0}/2$ and
$c_{k}=\alpha_{k}-\ii\beta_{k}$, for $k\in\{1,2,3,\ldots,\infty\}$. As was
shown in~\cite{CAoRFI}, the complex power series generated from these
coefficients,

\begin{equation}
  S(z)
  =
  \sum_{k=0}^{\infty}
  c_{k}z^{k},
\end{equation}

\noindent
always converges to an inner analytic function $w(z)$ within the open unit
disk,

\begin{equation}
  w(z)
  =
  u(\rho,\theta)+\ii v(\rho,\theta).
\end{equation}

\noindent
As was also shown in~\cite{CAoRFI}, the $\rho\to 1_{(-)}$ limit of the
real part $u(\rho,\theta)$ reproduces $f(\theta)$ at all points on the
unit circle where $w(z)$ does not have hard singularities. It does have
hard singularities at all points where $f(\theta)$ does, so we are led to
impose that these must be finite in number. However, in some special cases
$w(z)$ may have hard singularities at points where $f(\theta)$ does not,
and therefore we are led to assume independently that the number of hard
singularities of $w(z)$ is finite. For all integrable real functions
$f(\theta)$ that correspond to inner analytic functions $w(z)$ which have
at most a finite number of hard singularities on the unit circle, we have
that

\begin{equation}
  f(\theta)
  =
  \lim_{\rho\to 1_{(-)}}u(\rho,\theta),
\end{equation}

\noindent
almost everywhere. Since, being the real part of an analytic function, the
real function $u(\rho,\theta)$ is a harmonic function defined on the
plane, and thus satisfies the Laplace equation within the open unit disk,

\begin{equation}
  \nabla^{2}u(\rho,\theta)
  =
  0,
\end{equation}

\noindent
this construction establishes the existence of a solution of the Dirichlet
problem on the unit disk or, more precisely, the existence of a solution
of the Dirichlet boundary value problem of the Laplace equation on the
unit disk. Given the boundary condition $u(1,\theta)=f(\theta)$, the
solution is $u(\rho,\theta)$, which by construction satisfies the Laplace
equation within the open unit disk and which, also by construction,
assumes the values $f(\theta)$ on the unit circle, at least almost
everywhere.

Note that, since $f(\theta)$ may have isolated singular points where it
diverges to infinity, at which it is, therefore, not well defined, it is
clear that $u(\rho,\theta)$ can reproduce $f(\theta)$ only almost
everywhere. However, $u(\rho,\theta)$ may fail to reproduce $f(\theta)$ at
points other than its hard singularities, namely points where $w(z)$ has
hard singularities but $f(\theta)$ happens to have soft ones, due to the
way in which the complex singularities of $w(z)$ are oriented with respect
to the directions tangent to the unit circle at these singular points. In
this case the $\rho\to 1_{(-)}$ limit of $u(\rho,\theta)$ does not exist
at such points, and therefore at these points it is not possible to
recover the values of $f(\theta)$ in this way.

Note also that, if one introduces some removable singularities of
$f(\theta)$ at some points on the unit circle, then this does not change
the Fourier coefficients $\alpha_{0}$, $\alpha_{k}$ and $\beta_{k}$, for
$k\in\{1,2,3,\ldots,\infty\}$, since these are given by integrals, which
implies that it does not change the Taylor coefficients $c_{0}$ and
$c_{k}$, for $k\in\{1,2,3,\ldots,\infty\}$, and therefore that it also
does not change the corresponding inner analytic function $w(z)$. It
follows that $u(\rho,\theta)$ cannot reproduce $f(\theta)$ at these
points, if arbitrary real values are attributed to $f(\theta)$ at
them. Therefore, we are led to also assume that $f(\theta)$ has no such
removable singularities or, equivalently, we are led to assume that all
such removable singularities have been removed, and the function redefined
by continuity at these trivial singular points.

Here is, then, a complete and precise statement of the Dirichlet problem
on the unit disk, followed by the complete set of assumptions to be
imposed on $f(\theta)$ in order to ensure the existence of the solution of
that problem.

\begin{definition}\Colon\label{Def01}
  The Dirichlet Problem on the Unit Disk
\end{definition}

\noindent
Given the unit circle on the complex plane and a real function $f(\theta)$
defined on it, the existence problem of the Dirichlet boundary value
problem of the Laplace equation on the unit disk is to show that a
function $u(\rho,\theta)$ exists such that it satisfies

\begin{equation}
  \nabla^{2}u(\rho,\theta)
  =
  0,
\end{equation}

\noindent
within the open unit disk, and such that it also satisfies

\begin{equation}
  u(1,\theta)
  =
  f(\theta),
\end{equation}

\noindent
almost everywhere on the unit circle.

\vspace{2.6ex}

\noindent
In this section, using our results from previous papers, we will establish
the following theorem.

\begin{theorem}\Colon\label{Theo01}
  Given a real function $f(\theta)$ at the boundary of the unit disk, that
  satisfies the list of conditions described below, there is a solution
  $u(\rho,\theta)$ of the Dirichlet problem of the Laplace equation within
  the open unit disk, that assumes the values $f(\theta)$ almost
  everywhere at its boundary, the unit circle.
\end{theorem}

\begin{proof}\Colon
\end{proof}

\noindent
According to the construction introduced in~\cite{CAoRFI} and reviewed
above, which provides $u(\rho,\theta)$ starting from $f(\theta)$, the
function $u(\rho,\theta)$ that results from that construction is a
solution to this problem so long as $f(\theta)$ satisfies the following
set of conditions, which ensure that the construction of the inner
analytic function $w(z)$ from the real function $f(\theta)$ succeeds, and
that the real part $u(\rho,\theta)$ of $w(z)$ reproduces $f(\theta)$
almost everywhere oven the unit circle in the $\rho\to 1_{(-)}$ limit.
Apart from the global conditions that the real function $f(\theta)$ be a
Lebesgue-measurable function and that it have no removable singularities,
the conditions on $f(\theta)$ for this theorem are as follows.

\begin{enumerate}

\item The real function $f(\theta)$ is integrable on the unit circle.

\item The number of hard singularities of the corresponding inner analytic
  function $w(z)$ is finite.

\end{enumerate}

\noindent
This completes the proof of Theorem~\ref{Theo01}.

\vspace{2.6ex}

\noindent
Note that the last condition implies, in particular, that the number of
hard singularities of $f(\theta)$, where it is either discontinuous or
diverges to infinity, is also finite. Note also that, since the function
must be integrable, any hard singularities where it diverges to infinity
must be integrable ones, in the real sense of the terms involved. This
requires that these hard singularities be all isolated from each other, so
that there is a neighborhood around each one of them within which the two
lateral asymptotic limits of integrals can be considered. It is important
to emphasize that the conditions above over the real functions $f(\theta)$
include functions which are non-differentiable at any number of points,
discontinuous at a finite number of points, and unbounded at a finite
number of points, thus constituting a rather large set of boundary
conditions.

\section{Non-Integrable Real Functions on the Unit Circle}\label{Sec03}

In a previous paper~\cite{CAoRFIV} we showed that the correspondence
between real functions and inner analytic functions established
in~\cite{CAoRFI} can be extended to non-integrable real functions, so long
as these functions are locally integrable almost everywhere, and so long
as the non-integrable hard singularities of the corresponding inner
analytic functions have finite degrees of hardness. The definition of
local integrability almost everywhere on the unit circle is as follows.

\begin{definition}\Colon\label{Def02}
  A real function $f(\theta)$ is locally integrable almost everywhere on
  the unit circle if it is integrable on every closed interval
  $[\theta_{\ominus},\theta_{\oplus}]$ contained within that domain, that
  does not contain any of the points where the function has non-integrable
  hard singularities, of which there is a finite number.
\end{definition}

\noindent
Although the construction used in this case, which is given
in~\cite{CAoRFIV}, is considerably more involved than the one for the case
of integrable real functions, it is still true that given such a
non-integrable real function one can define a unique inner analytic
function $w(z)$ that corresponds to it, as well as a unique and complete
set of complex Taylor coefficients $c_{0}=\alpha_{0}/2$ and
$c_{k}=\alpha_{k}-\ii\beta_{k}$, for $k\in\{1,2,3,\ldots,\infty\}$, and
thus a corresponding unique and complete set of Fourier coefficients
$\alpha_{0}$, $\alpha_{k}$ and $\beta_{k}$, for
$k\in\{1,2,3,\ldots,\infty\}$, that are associated to it, despite the fact
that the real function is not integrable. From the real part of this inner
analytic function one can, once again, recover the real function almost
everywhere by taking the $\rho\to 1_{(-)}$ limit to the unit circle.
Therefore we have at hand all that we need in order to implement the proof
of existence in this more general case.

\vspace{2.6ex}

\noindent
In this section, using again our results from previous papers, we will
establish the following theorem.

\begin{theorem}\Colon\label{Theo02}
  Given a real function $f(\theta)$ at the boundary of the unit disk, that
  satisfies the list of conditions described below, there is a solution
  $u(\rho,\theta)$ of the Dirichlet problem of the Laplace equation within
  the open unit disk, that assumes the values $f(\theta)$ almost
  everywhere at its boundary, the unit circle.
\end{theorem}

\begin{proof}\Colon
\end{proof}

\noindent
The argument is the same as the one used before for Theorem~\ref{Theo01}
in Section~\ref{Sec02}, in the case of integrable real functions, but
using now the construction presented in~\cite{CAoRFIV}, instead of the one
presented in~\cite{CAoRFI}. Due to this, the only change with respect to
that previous case is that our list of conditions on the real functions
can now be upgraded to the following, still including the previous case.

\begin{enumerate}

\item The real function $f(\theta)$ is locally integrable almost
  everywhere on the unit circle, including the cases in which this
  function is globally integrable there.

\item The number of hard singularities of the corresponding inner analytic
  function $w(z)$ is finite.

\item The hard singularities of the corresponding inner analytic function
  $w(z)$ have finite degrees of hardness.

\end{enumerate}

\noindent
This completes the proof of Theorem~\ref{Theo02}.

\vspace{2.6ex}

\noindent
In this way we have generalized the proof of existence of the Dirichlet
problem from boundary conditions given by integrable real functions to
others given by a certain class of non-integrable real functions. Note
that the degree of hardness in the previous case, that of borderline hard
singularities of integrable real functions, is simply zero. In other
words, if all the hard singularities are borderline hard ones, then the
function is simply integrable. Therefore, this theorem is a strict
generalization of the previous one, and contains it.

\section{Conformal Transformations to Other Curves}\label{Sec04}

As we will see in the subsequent sections, it is possible to extend the
proof of existence of the Dirichlet problem to boundaries other than the
unit circle, through the use of conformal transformations. Therefore, as a
preliminary to the proof of further versions of the existence theorem, in
this section we will describe such conformal transformations and examine
some of their well-known properties, targeting their use here. In order to
do this, consider two complex variables $z_{a}$ and $z_{b}$ and the
corresponding complex planes, a complex analytic function $\gamma(z)$
defined on the complex plane $z_{a}$ with values on the complex plane
$z_{b}$, and its inverse function, which is a complex analytic function
$\gamma^{(-1)}(z)$ defined on the complex plane $z_{b}$ with values on the
complex plane $z_{a}$,

\noindent
\begin{eqnarray}
  z_{b}
  & = &
  \gamma(z_{a}),
  \nonumber\\
  z_{a}
  & = &
  \gamma^{(-1)}(z_{b}).
\end{eqnarray}

\noindent
Consider a bounded and simply connected open region $S_{a}$ on the complex
plane $z_{a}$ and its image $S_{b}$ under $\gamma(z)$, which is a similar
region on the complex plane $z_{b}$. It can be shown that if $\gamma(z)$
is analytic on $S_{a}$, is invertible there, and its derivative has no
zeros there, then its inverse function $\gamma^{(-1)}(z)$ has these same
three properties on $S_{b}$, and the mapping between the two complex
planes established by $\gamma(z)$ and $\gamma^{(-1)}(z)$ is conformal, in
the sense that it preserves the angles between oriented curves at points
where they cross each other. Note that this mapping is a bijection between
the two regions, and establishes an equivalence relation that can be
extended in a transitive way to other regions.

Consider now that the regions under consideration are the interiors of
simple closed curves. One of these curves will be the unit circle $C_{a}$
on the complex plane $z_{a}$, and the other will be a given curve $C_{b}$
on the complex plane $z_{b}$. We will assume that the curve $C_{b}$ has
finite total length, for two reasons, one being to ensure that the
interior of the curve is a bounded set, and the other being to ensure that
the integrals of real functions over the curve $C_{b}$ are integrals over
a finite-length, compact domain. Since $\gamma(z_{a})$, being analytic, is
in particular a continuous function, the image on the $z_{b}$ plane of the
unit circle $C_{a}$ on the $z_{a}$ plane must be a continuous closed curve
$C_{b}$. We can also see that $C_{b}$ must be a simple curve, because the
fact that $\gamma(z_{a})$ is invertible on $C_{a}$ means that it cannot
have the same value at two different points of $C_{a}$, and therefore no
two points of $C_{b}$ can be the same. Consequently, the curve $C_{b}$
cannot self-intersect.

We thus see that, so far, we are restricted to simple closed curves
$C_{b}$ with finite total lengths. However, there are further limitations
on the curves, implied by our hypotheses. Since the transformation is
conformal and thus preserves angles, it follows that in this case the
smooth unit circle $C_{a}$ is mapped onto another equally differentiable
circuit $C_{b}$. One can see this by considering the angles between
tangents to the curve $C_{b}$ at pairs of neighboring points, the
corresponding elements on the curve $C_{a}$, and the limit of these angles
when the two points tend to each other, given that the transformation is
conformal. Therefore, with such limitations one cannot map the unit circle
onto a square or any other polygon. This limitation can be lifted by
allowing the derivative of $\gamma(z_{a})$ to have a finite number of
isolated zeros on the curve $C_{a}$, which then implies that the
derivative of $\gamma^{(-1)}(z_{b})$ will have a finite number of
corresponding isolated singular points on $C_{b}$.

Let us assume that the unit circle $C_{a}$ is described by the real
arc-length parameter $\theta$ on the $z_{a}$ plane, and that the curve
$C_{b}$ is described by a corresponding real parameter $\lambda$ on the
$z_{b}$ plane. Let us assume also that $\lambda$ is chosen in such a way
that $|d\lambda|=|dz_{b}|$ over the curve $C_{b}$, just as
$|d\theta|=|dz_{a}|$ over $C_{a}$, which means that $\lambda$ is also an
arc-length parameter. Since every point $z_{a}$ on the curve $C_{a}$ is
mapped by the conformal transformation onto a corresponding point $z_{b}$
on the curve $C_{b}$, and since a point $z_{a}$ on $C_{a}$ is described by
a certain value of $\theta$, while the corresponding point $z_{b}$ on
$C_{b}$ is described by a certain value of $\lambda$, it is clear that the
complex conformal transformation induces a corresponding real
transformation between the values of $\theta$ and the values of $\lambda$,

\noindent
\begin{eqnarray}
  z_{b}
  & = &
  \gamma(z_{a})
  \;\;\;\Rightarrow
  \nonumber\\
  \lambda
  & = &
  g(\theta),
\end{eqnarray}

\noindent
where the real function $g(\theta)$ is continuous, differentiable and
invertible on $C_{a}$. We will refer to the function $g(\theta)$ as the
real transformation induced on the curve $C_{a}$ by the complex conformal
transformation $\gamma(z_{a})$. The same is true for the inverse
transformation, which induces the inverse function of $g(\theta)$, on the
curve $C_{b}$,

\noindent
\begin{eqnarray}
  z_{a}
  & = &
  \gamma^{(-1)}(z_{b})
  \;\;\;\Rightarrow
  \nonumber\\
  \theta
  & = &
  g^{(-1)}(\lambda),
\end{eqnarray}

\noindent
where the real function $g^{(-1)}(\lambda)$ is continuous and invertible
on $C_{b}$, and also differentiable so long as $C_{b}$ is a differentiable
curve. Before we proceed, we must now consider in more detail the question
of what is the set of curves $C_{b}$ for which the structure described
above can be set up. We assume that this curve is a simple closed curve of
finite total length, and the question is whether or not this structure can
be set up for an arbitrary such curve. Given the curve $C_{b}$, the only
additional objects we need in order to do this is the conformal mapping
$\gamma(z_{a})$ and its inverse $\gamma^{(-1)}(z_{b})$, between that curve
and the unit circle $C_{a}$.

The existence of these transformation functions can be ensured as a
consequence of the famous Riemann mapping theorem, and of the associated
results relating to conformal mappings between regions of the complex
plane. According to that theorem, a conformal transformation such as the
one we described here exists between any bounded simply connected open set
of the plane and the open unit disk. In addition to this, one can show
that this conformal mapping can be extended to the respective boundaries
as a continuous function so long as the boundary curve $C_{b}$ satisfies a
certain condition~\cite{RMPQiu}.

The condition on $C_{b}$ that implies the existence of the continuous
extension to the boundary is that every point on that curve be what is
called in the relevant literature a {\em simple point}. This means that no
point of $C_{b}$ can be a {\em multiple point}, which in essence is a
point on the boundary that is accessible from the interior via two or more
independent continuous paths contained in the interior, that cannot be
continuously deformed into each other without crossing the boundary. We
can see, therefore, that this condition has a topological character. Note
that the presence of a multiple point on the boundary means that, even if
the open set under consideration is simply connected, its {\em closure}
will {\em not} be. Therefore, one way to formulate this condition is to
simply state that the closure of the bounded simply connected open set
must also be simply connected.

Since the existence of a multiple point at the boundary means that this
boundary is not a simple curve, it follows that, under the limitations
over $C_{b}$ that we have here, the conformal mapping on the open unit
disk can always be continuously extended to the unit circle, and hence
from the interior of the curve $C_{b}$ to that curve, which is mapped from
the unit circle. Therefore, we conclude that it is a known fact that such
a conformal transformation exists for all possible simple closed curves
$C_{b}$ with finite total lengths, and in particular for all such curves
which are also differentiable, in which case the extension is also
differentiable on the unit circle $C_{a}$. It is therefore not necessary
to impose explicitly any additional hypotheses about the existence of the
conformal transformation, regardless of whether or not the curves $C_{b}$
under consideration are differentiable.

Let us close this section with a discussion of the nature of the
singularities that appear in the case of simple closed curves $C_{b}$
which are not differentiable at a finite set of points $z_{b,i}$, for
$i\in\{1,\ldots,N\}$. The additional difficulty that appears in this case
stems from the fact that, if the curve $C_{b}$ is not differentiable at
the points $z_{b,i}$, then the derivative of the transformation
$\gamma(z_{a})$ has isolated zeros at the corresponding points $z_{a,i}$
on the curve $C_{a}$, and therefore the derivative of the inverse
transformation $\gamma^{(-1)}(z_{b})$ has isolated hard singularities at
the points $z_{b,i}$. In order to see how this comes about we start by
noting that, since we have that

\begin{equation}
  \gamma^{(-1)}(\gamma(z_{a}))
  =
  z_{a},
\end{equation}

\noindent
for all $z_{a}$ on the closed unit disk, differentiating this equation we
get, due to the chain rule,

\begin{equation}
  \frac{d\gamma^{(-1)}}{dz_{b}}\,
  \frac{d\gamma}{dz_{a}}\,
  =
  1.
\end{equation}

\noindent
Therefore, to every point $z_{a,i}$ on $C_{a}$ where the derivative of the
transformation has a zero corresponds a point $z_{b,i}$ on $C_{b}$ where
the derivative of the inverse transformation diverges to infinity. Taking
absolute values we have, in terms of the arc-length parameters $\theta$
and $\lambda$,

\noindent
\begin{eqnarray}
  \left|
    \frac{d\gamma^{(-1)}}{dz_{b}}
  \right|
  \left|
    \frac{d\gamma}{dz_{a}}
  \right|
  & = &
  \left|
    \frac{dz_{a}}{dz_{b}}
  \right|
  \left|
    \frac{dz_{b}}{dz_{a}}
  \right|
  \nonumber\\
  & = &
  \left|
    \frac{d\theta}{d\lambda}
  \right|
  \left|
    \frac{d\lambda}{d\theta}
  \right|,
\end{eqnarray}

\noindent
which implies that

\begin{equation}\label{EQInvDir}
  \left|
    \frac{d\theta}{d\lambda}
  \right|
  \left|
    \frac{d\lambda}{d\theta}
  \right|
  =
  1.
\end{equation}

\noindent
In fact, the real function in the left-hand side of this last equation has
a removable singularity at every point where the first derivative in the
product diverges and the second one is zero. Consequently, they can be
removed by simply redefining the product by continuity at these points.
When we approach one of the points $z_{a,i}$, which are characterized by
the values $\theta_{i}$ of the parameter $\theta$, along the curve
$C_{a}$, we have that

\noindent
\begin{eqnarray}
  \theta
  & \to &
  \theta_{i},
  \nonumber\\
  \lambda
  & \to &
  \lambda_{i}
  \nonumber\\
  & \Downarrow &
  \nonumber\\
  \left|
    \frac{d\lambda}{d\theta}
  \right|
  & \to &
  0,
  \nonumber\\
  \left|
    \frac{d\theta}{d\lambda}
  \right|
  & \to &
  \infty,
\end{eqnarray}

\noindent
where the corresponding points $z_{b,i}$ are characterized by the values
$\lambda_{i}$ of the parameter $\lambda$, along the curve $C_{b}$.
Although the derivative $d\theta/d\lambda$ does, therefore, have hard
singularities at $z_{b,i}$, we can show that these are still integrable
singularities. We simply integrate the expressions in either side of
Equation~(\ref{EQInvDir}) absolutely over $C_{a}$, thus obtaining

\begin{equation}
  \oint_{C_{a}}|d\theta|\,
  \left|
    \frac{d\theta}{d\lambda}
  \right|
  \left|
    \frac{d\lambda}{d\theta}
  \right|
  =
  2\pi.
\end{equation}

\noindent
If we now change variables in this integral from $\theta$ to $\lambda$, we
get the integral over $C_{b}$

\begin{equation}
  \oint_{C_{b}}|d\lambda|\,
  \left|
    \frac{d\theta}{d\lambda}
  \right|
  =
  2\pi.
\end{equation}

\noindent
This shows that the real function appearing as the integrand in this
integral is an integrable real function on $C_{b}$. Therefore the hard
singularities where the derivative $d\theta/d\lambda$ diverges to infinity
are integrable hard singularities, which therefore have degree of hardness
zero. These are also referred to as borderline hard singularities. Note
that, as a consequence, the corresponding singularities of the inverse
real transformation $g^{(-1)}(\lambda)$ itself must be soft ones, with
degrees of softness equal to one.

\section{Integrable Real Functions on Differentiable Curves}\label{Sec05}

In this section we will show how one can generalize the proof of existence
of the solution of the Dirichlet problem on the unit disk, given by
Theorem~\ref{Theo01} in Section~\ref{Sec02}, to the case in which we have,
as the boundary condition, integrable real functions defined on boundaries
given by differentiable curves on the plane. In order to do this, the
first thing we must do here is to establish the precise definition of the
Dirichlet problem in this case.

\begin{definition}\Colon\label{Def03}
  The Dirichlet Problem on a Given Curve and its Interior
\end{definition}

\noindent
Given a simple closed curve $C$ on the complex plane, described by a real
arc-length parameter $\lambda$, and a real function $f(\lambda)$ defined
on it, the existence problem of the Dirichlet boundary value problem of
the Laplace equation on this curve and its interior is to show that a
function $u(x,y)$ exists such that it satisfies

\begin{equation}
  \nabla^{2}u(x,y)
  =
  0,
\end{equation}

\noindent
within the interior of $C$, and such that it also satisfies

\begin{equation}
  u(x,y)
  =
  f(\lambda),
\end{equation}

\noindent
for $z=x+\ii y$ on $C$, thus corresponding to $\lambda$, almost everywhere
over that curve.

\vspace{2.6ex}

\noindent
Therefore, using the notation established in Section~\ref{Sec04}, let
$C_{b}$ be a differentiable simple closed curve on the complex $z_{b}$
plane, with finite total length, and let us assume that a Dirichlet
boundary value problem for the Laplace equation is given on the region
whose boundary is the curve $C_{b}$, that is, let there be given also an
integrable real function $f_{b}(\lambda)$ on $C_{b}$, that is, a function
such that the integral

\begin{equation}
  \oint_{C_{b}}d\lambda\,
  f_{b}(\lambda)
\end{equation}

\noindent
exists and is finite. The problem is then to establish the existence of a
function $u_{b}(x,y)$ that satisfies $\nabla^{2}u_{b}(x,y)=0$ in the
interior of the curve $C_{b}$ and that assumes the values $f_{b}(\lambda)$
almost everywhere over that curve. In order to do this using the conformal
transformation $\gamma(z)$ from the complex plane $z_{a}$ to the complex
plane $z_{b}$, we start by constructing a corresponding Dirichlet problem
on the unit disk of the $z_{a}$ plane, using the mapping between $z_{b}$
and $z_{a}$ provided by the conformal transformation $\gamma(z_{a})$ and
its inverse $\gamma^{(-1)}(z_{b})$. We define a corresponding real
function $f_{a}(\theta)$ on the unit circle $C_{a}$ by simply transferring
the values of $f_{b}(\lambda)$ through the use of the conformal mapping
from point to point,

\noindent
\begin{eqnarray}
  f_{a}(\theta)
  & = &
  f_{b}(\lambda)
  \nonumber\\
  & = &
  f_{b}(g(\theta)),
\end{eqnarray}

\noindent
where $\theta$ describes a point on $C_{a}$ given by the complex number
$z_{a}$, $\lambda$ describes the corresponding point on $C_{b}$ given by
the complex number $z_{b}=\gamma(z_{a})$, and $g(\theta)$ is the induced
real transformation, so that we have that $\lambda=g(\theta)$ and that
$\theta=g^{(-1)}(\lambda)$. We will start by establishing the following
preliminary fact about the real function $f_{a}(\theta)$ defined as above
from an integrable real function $f_{b}(\lambda)$.

\begin{lemma}\Colon\label{Lemma01}
  Given a real function $f_{b}(\lambda)$ which is integrable on the
  differentiable simple closed curve $C_{b}$ of finite total length, it
  follows that the corresponding function $f_{a}(\theta)$ defined on the
  unit circle $C_{a}$ by $f_{a}(\theta)=f_{b}(\lambda)$ is integrable on
  that circle.
\end{lemma}

\noindent
Given that $f_{b}(\lambda)$ is integrable on $C_{b}$, we must show that
$f_{a}(\theta)$ defined by the composition of $f_{b}(\lambda)$ with
$g(\theta)$, is also integrable, that is, it is integrable on $C_{a}$.
Since all real functions under discussion here are assumed to be
Lebesgue-measurable, and since for such measurable functions defined on
compact domains integrability and absolute integrability are equivalent
conditions~\cite{RARudin,RARoyden}, given that $f_{b}(\lambda)$ is
integrable on $C_{b}$ we have that the integral

\begin{equation}
  \oint_{C_{b}}|d\lambda|\,
  |f_{b}(\lambda)|
\end{equation}

\noindent
exists and is finite. We must now show that the integral

\begin{equation}\label{EQAbsInt1}
  \oint_{C_{a}}|d\theta|\,
  |f_{a}(\theta)|
\end{equation}

\noindent
exists and is finite, which is equivalent to the statement that
$f_{a}(\theta)$ is integrable on $C_{a}$. Changing variables on this
integral from $\theta$ to $\lambda$, and using the fact that by definition
we have that $f_{a}(\theta)=f_{b}(\lambda)$, we obtain

\begin{equation}\label{EQAbsInt2}
  \oint_{C_{a}}|d\theta|\,
  |f_{a}(\theta)|
  =
  \oint_{C_{b}}|d\lambda|\,
  \left|
    \frac{d\theta}{d\lambda}
  \right|
  |f_{b}(\lambda)|.
\end{equation}

\noindent
Since $f_{b}(\lambda)$ is integrable on $C_{b}$, and since the absolute
value of the derivative shown exists and is finite, given that we have

\noindent
\begin{eqnarray}
  \left|
    \frac{d\theta}{d\lambda}
  \right|
  & = &
  \left|
    \frac{dz_{a}}{dz_{b}}
  \right|
  \nonumber\\
  & = &
  \left|
    \frac{d\gamma^{(-1)}(z_{b})}{dz_{b}}
  \right|,
\end{eqnarray}

\noindent
where $\gamma^{(-1)}(z_{b})$ is analytic on $C_{b}$ and therefore
differentiable there, it follows that the absolute value of the derivative
which appears in the integrand on the right-hand side of
Equation~(\ref{EQAbsInt2}) is a limited real function on $C_{b}$. Since
$f_{b}(\lambda)$ is integrable on $C_{b}$, from this it follows that the
whole integrand of the integral on the right-hand side of
Equation~(\ref{EQAbsInt2}), which is the product of a limited real
function with an integrable real function, is itself an integrable real
function on $C_{b}$, so that we may conclude that the integral in
Equation~(\ref{EQAbsInt1}) exists and is finite, and therefore that
$f_{a}(\theta)$ is an integrable real function on $C_{a}$. This
establishes Lemma~\ref{Lemma01}.

\vspace{2.6ex}

\noindent
In this section, using the results from the previous sections, we will
establish the following theorem.

\begin{theorem}\Colon\label{Theo03}
  Given a differentiable simple closed curve $C$ of finite total length on
  the complex plane $z=x+\ii y$, given the invertible conformal
  transformation $\gamma(z)$ whose derivative has no zeros on the closed
  unit disk, that maps it from the unit circle, and given a real function
  $f(\lambda)$ on that curve, that satisfies the list of conditions
  described below, there is a solution $u(x,y)$ of the Dirichlet problem
  of the Laplace equation within the interior of that curve, that assumes
  the given values $f(\lambda)$ almost everywhere on the curve.
\end{theorem}

\begin{proof}\Colon
\end{proof}

\noindent
The proof consists of using the conformal transformation between the
closed unit disk and the union of the curve $C$ with its interior, a
transformation which according to the analysis in Section~\ref{Sec04}
always exists, to map the given boundary condition on $C$ onto a
corresponding boundary condition on the unit circle, then using the proof
of existence established before by Theorem~\ref{Theo01} in
Section~\ref{Sec02} for the closed unit disk to establish the existence of
the solution of the corresponding Dirichlet problem on that disk, and
finally using once more the conformal transformation to map the resulting
solution back to $C$ and its interior, showing in the process that one
obtains in this way the solution of the Dirichlet problem there. The list
of conditions on the real functions is now the following.

\begin{enumerate}

\item The real function $f(\lambda)$ is integrable on $C$.

\item The number of hard singularities on the unit circle of the
  corresponding inner analytic function $w(z)$ on the unit disk is finite.

\end{enumerate}

\noindent
According to the preliminary result established in Lemma~\ref{Lemma01}, if
the real function $f_{b}(\lambda)$ satisfies these conditions on $C_{b}$,
then $f_{a}(\theta)$ is an integrable real function on $C_{a}$. Therefore,
due to the existence theorem of the Dirichlet problem on the unit disk of
the plane $z_{a}$, which was established by Theorem~\ref{Theo01} in
Section~\ref{Sec02}, we know that there is an inner analytic function
$w_{a}(z_{a})$ such that its real part $u_{a}(\rho,\theta)$ is harmonic
within the open unit disk and also satisfies
$u_{a}(1,\theta)=f_{a}(\theta)$ almost everywhere at the boundary $C_{a}$.
Now, by composing $w_{a}(z_{a})$ with the inverse conformal transformation
$\gamma^{(-1)}(z_{b})$ we get on the $z_{b}$ plane the complex function

\noindent
\begin{eqnarray}
  w_{b}(z_{b})
  & = &
  w_{a}(z_{a})
  \nonumber\\
  & = &
  w_{a}\!\left(\gamma^{(-1)}(z_{b})\right),
\end{eqnarray}

\noindent
which corresponds to simply transferring back the values of
$w_{a}(z_{a})$, by the use of the conformal mapping from point to point,
while we also have, of course, the corresponding inverse real
transformation at the boundary,

\noindent
\begin{eqnarray}
  f_{b}(\lambda)
  & = &
  f_{a}(\theta)
  \nonumber\\
  & = &
  f_{a}\!\left(g^{(-1)}(\lambda)\right).
\end{eqnarray}

\noindent
Given that the composition of two analytic functions is also analytic, in
their chained domain of analyticity, and since $\gamma^{(-1)}(z_{b})$ is
analytic in the interior of the curve $C_{b}$, and also since
$w_{a}(z_{a})$ is analytic in the interior of the curve $C_{a}$, we
conclude that $w_{b}(z_{b})$ is analytic in the interior of the curve
$C_{b}$. Therefore, the real part $u_{b}(x,y)$ of $w_{b}(z_{b})$ is
harmonic and thus satisfies

\begin{equation}
  \nabla^{2}u_{b}(x,y)
  =
  0
\end{equation}

\noindent
in the interior of $C_{b}$, while by construction the fact that we have
$u_{a}(\rho,\theta)=f_{a}(\theta)$ for the points given by
$z_{a}=\rho\exp(\ii\theta)$ almost everywhere on $C_{a}$, where $\rho=1$,
implies that we also have

\begin{equation}
  u_{b}(x,y)=f_{b}(\lambda),
\end{equation}

\noindent
for the corresponding points given by $z_{b}=x+\ii y$ almost everywhere on
$C_{b}$, and which thus correspond to $\lambda$. This establishes the
existence, by construction, of the solution of the Dirichlet problem on
the $z_{b}$ plane, under our current hypotheses. This completes the proof
of Theorem~\ref{Theo03}.

\vspace{2.6ex}

\noindent
In this way we have generalized the proof of existence of the Dirichlet
problem from the unit circle to all differentiable simple closed curves
with finite total lengths on the plane, for boundary conditions given by
integrable real functions.

\section{Integrable Real Functions on Non-Differentiable
  Curves}\label{Sec06}

In this section we will extend the existence theorem of the Dirichlet
problem on the unit disk, given by Theorem~\ref{Theo01} in
Section~\ref{Sec02}, to regions bounded by simple closed curves $C_{b}$
which are not differentiable at a finite set of points $z_{b,i}$, for
$i\in\{1,\ldots,N\}$. We will still use the conformal transformation known
to exist between the open unit disk and the interior of any such curve, as
well as its continuous extension to the respective boundaries, where the
extension is also differentiable almost everywhere, with the exception of
a finite set of singularities of the inverse conformal transformation, at
the points $z_{b,i}$, where the inverse conformal transformation still
exists but is not differentiable. Note that in Section~\ref{Sec04} we
established the existence of the conformal mapping $\gamma(z)$ for all
simple closed curves $C_{b}$ with finite total lengths, regardless of
whether or not the curves are differentiable.

The additional difficulty that appears in this case stems from the fact
that, if the curve $C_{b}$ is not differentiable at the points $z_{b,i}$,
then the derivative of the transformation $\gamma(z_{a})$ has isolated
zeros at the corresponding points $z_{a,i}$ on the curve $C_{a}$, and
therefore the derivative of the inverse transformation
$\gamma^{(-1)}(z_{b})$ has isolated hard singularities at the points
$z_{b,i}$, as was discussed in Section~\ref{Sec04}. This will require that
we impose one additional limitation on the real functions giving the
boundary conditions, namely that any integrable hard singularities where
they diverge to infinity do not coincide with any of the points $z_{b,i}$.

The precise definition of the Dirichlet problem in this case is the same
one given in Definition~\ref{Def03}, in Section~\ref{Sec05}. We will start
by establishing the following preliminary fact about the real function
$f_{a}(\theta)$ defined from an integrable real function $f_{b}(\lambda)$.

\begin{lemma}\Colon\label{Lemma02}
  Given a real function $f_{b}(\lambda)$ which is integrable on the simple
  closed curve $C_{b}$ of finite total length, which is not differentiable
  at a finite set of points $z_{b,i}$, for $i\in\{1,\ldots,N\}$, and given
  that the integrable hard singularities of $f(\lambda)$ where it diverges
  to infinity are not located at any of the points $z_{b,i}$ where the
  curve is non-differentiable, it follows that the corresponding function
  $f_{a}(\theta)$ defined on the unit circle $C_{a}$ by
  $f_{a}(\theta)=f_{b}(\lambda)$ is integrable on that circle.
\end{lemma}

\noindent
Given that $f_{b}(\lambda)$ is integrable on $C_{b}$, we must decide
whether or not $f_{a}(\theta)$ is also integrable, that is, whether it is
integrable on $C_{a}$. Once again, given that $f_{b}(\lambda)$ is
integrable on $C_{b}$ we have that the integral

\begin{equation}
  \oint_{C_{b}}|d\lambda|\,
  |f_{b}(\lambda)|
\end{equation}

\noindent
exists and is finite. We must now determine whether or not the integral

\begin{equation}
  \oint_{C_{a}}|d\theta|\,
  |f_{a}(\theta)|
\end{equation}

\noindent
exists and is finite, which is equivalent to the statement that
$f_{a}(\theta)$ is integrable on $C_{a}$. Changing variables on this
integral from $\theta$ to $\lambda$ we obtain once again

\begin{equation}
  \oint_{C_{a}}|d\theta|\,
  |f_{a}(\theta)|
  =
  \oint_{C_{b}}|d\lambda|\,
  \left|
    \frac{d\theta}{d\lambda}
  \right|
  |f_{b}(\lambda)|.
\end{equation}

\noindent
Since both the absolute value of the derivative shown and the function
$f_{b}(\lambda)$ are integrable on $C_{b}$, and since the integrable
borderline hard singular points where either one of these two real
functions diverges to infinity do not coincide, we have that the integrand
of the integral on the right-hand side of this equation is also an
integrable real function, and thus that the integral exists and is
finite. We can see that the integrand is an integrable real function
because around each integrable hard singular point of either one of the
two real functions involved there is a neighborhood where the other real
function is limited. Since the product of a limited real function with an
integrable real function is also an integrable real function, we may
conclude that the integrand is locally integrable {\em everywhere} on
$C_{b}$, and therefore globally integrable there, so that the integral
above exists and is finite. It thus follows that $f_{a}(\theta)$ is an
integrable real function on $C_{a}$. This establishes Lemma~\ref{Lemma02}.

\vspace{2.6ex}

\noindent
In this section, using again the results from the previous sections, we
will establish the following theorem.

\begin{theorem}\Colon\label{Theo04}
  Given a simple closed curve $C$ of finite total length on the complex
  plane $z=x+\ii y$, which is non-differentiable at a given finite set of
  points $z_{i}$, for $i\in\{1,\ldots,N\}$, given the conformal
  transformation $\gamma(z)$ that maps it from the unit circle, whose
  derivative has zeros on the unit circle at the corresponding set of
  points, and given a real function $f(\lambda)$ on that curve, that
  satisfies the list of conditions described below, there is a solution
  $u(x,y)$ of the Dirichlet problem of the Laplace equation within the
  interior of that curve, that assumes the given values $f(\lambda)$
  almost everywhere on the curve.
\end{theorem}

\newpage

\begin{proof}\Colon
\end{proof}

\noindent
Just as in the previous case, in Section~\ref{Sec05}, the proof consists
of using the conformal transformation between the closed unit disk and the
union of the curve $C$ with its interior, which according to the analysis
in Section~\ref{Sec04} always exists, to map the given boundary condition
on $C$ onto a corresponding boundary condition on the unit circle, then
using the proof of existence established before by Theorem~\ref{Theo01} in
Section~\ref{Sec02} for the closed unit disk to establish the existence of
the solution of the corresponding Dirichlet problem on that disk, and
finally using once more the conformal transformation to map the resulting
solution back to $C$ and its interior, thus obtaining the solution of the
original Dirichlet problem. The list of conditions on the real functions
is now the following.

\begin{enumerate}

\item The real function $f(\lambda)$ is integrable on $C$.

\item The number of hard singularities on the unit circle of the
  corresponding inner analytic function $w(z)$ on the unit disk is finite.

\item The integrable hard singularities of $f(\lambda)$ where it diverges
  to infinity are not located at any of the points where the curve $C$ is
  non-differentiable.

\end{enumerate}

\noindent
The rest of the proof is identical to that of the previous case, in
Section~\ref{Sec05}. Therefore, once again we may conclude that, due to
the existence theorem of the Dirichlet problem on the unit disk of the
plane $z_{a}$, which was established by Theorem~\ref{Theo01} in
Section~\ref{Sec02}, we know that there is an inner analytic function
$w_{a}(z_{a})$ such that its real part $u_{a}(\rho,\theta)$ is harmonic
within the open unit disk and satisfies $u_{a}(1,\theta)=f_{a}(\theta)$
almost everywhere at the boundary $C_{a}$. Just as in Section~\ref{Sec05},
we get on the $z_{b}$ plane the complex function $w_{b}(z_{b})$ which is
analytic in the interior of the curve $C_{b}$. Therefore, the real part
$u_{b}(x,y)$ of $w_{b}(z_{b})$ is harmonic and thus satisfies

\begin{equation}
  \nabla^{2}u_{b}(x,y)
  =
  0
\end{equation}

\noindent
in the interior of $C_{b}$, while we also have that

\begin{equation}
  u_{b}(x,y)=f_{b}(\lambda),
\end{equation}

\noindent
almost everywhere on $C_{b}$. This establishes the existence, by
construction, of the solution of the Dirichlet problem on the $z_{b}$
plane, under our current hypotheses. This completes the proof of
Theorem~\ref{Theo04}.

\vspace{2.6ex}

\noindent
In this way we have generalized the proof of existence of the Dirichlet
problem from the unit circle to all simple closed curves with finite total
lengths on the plane, that can be either differentiable or
non-differentiable on at most a finite set of points, still for boundary
conditions given by integrable real functions.

\section{Non-Integrable Real Functions on Differentiable
  Curves}\label{Sec07}

In this section we will show how one can generalize the proof of existence
of the Dirichlet problem on the unit disk, given by Theorem~\ref{Theo02}
in Section~\ref{Sec03}, to the case in which we have, as the boundary
condition, non-integrable real functions $f_{b}(\lambda)$ defined on
boundaries given by differentiable curves $C_{b}$ on the plane. We will be
able to do this if the non-integrable real functions, despite being
non-integrable over the whole curves $C_{b}$, are however locally
integrable almost everywhere on those curves, and if, in addition to this,
the non-integrable hard singularities of the inner analytic functions
involved have finite degrees of hardness. The definition of the concept of
local integrability almost everywhere is similar to that given for the
unit circle by Definition~\ref{Def02}, in Section~\ref{Sec03}. In our case
here the precise definition of local integrability almost everywhere is as
follows.

\begin{definition}\Colon\label{Def04}
  A real function $f(\lambda)$ is locally integrable almost everywhere on
  the curve $C$ described by the arc-length parameter $\lambda$ if it is
  integrable on every closed interval
  $[\lambda_{\ominus},\lambda_{\oplus}]$ contained within that domain,
  that does not contain any of the points where the function has
  non-integrable hard singularities, of which there is a finite number.
\end{definition}

\noindent
The proof will follow the general lines of the one given for integrable
real functions in Section~\ref{Sec05}, with the difference that, since the
real functions $f_{b}(\lambda)$ are assumed to be non-integrable on
$C_{b}$, but locally integrable almost everywhere there, instead of
showing that the corresponding functions $f_{a}(\theta)$ on the unit
circle $C_{a}$ are integrable there, we will show that they are locally
integrable almost everywhere there. In addition to this, instead of using
the result for integrable real function on the unit circle, which was
given by Theorem~\ref{Theo01} in Section~\ref{Sec02}, we will use the
corresponding result for non-integrable real functions which are locally
integrable almost everywhere on the unit circle, which was given by
Theorem~\ref{Theo02} in Section~\ref{Sec03}. Since that result depends
also on the non-integrable hard singularities of the real functions having
finite degrees of hardness, we will also show that the hypothesis that the
functions $f_{b}(\lambda)$ have this property implies that the
corresponding functions $f_{a}(\theta)$ have the same property as well. In
order to do this we will use the technique of {\em piecewise integration}
which was introduced and employed in~\cite{CAoRFIV}, where it played a
crucial role.

We will start by showing the following preliminary fact about a real
function $f_{a}(\theta)$ defined from a real function $f_{b}(\lambda)$
which is locally integrable almost everywhere on $C_{b}$, and which has
non-integrable hard singularities at the finite set of points $z_{b,j}$,
for $j\in\{1,\ldots,M\}$.

\begin{lemma}\Colon\label{Lemma03}
  Given a real function $f_{b}(\lambda)$ which is integrable on a given
  closed interval $I_{b}$ on $C_{b}$, it follows that the corresponding
  function $f_{a}(\theta)$ defined on the unit circle $C_{a}$ by
  $f_{a}(\theta)=f_{b}(\lambda)$ is integrable on the corresponding closed
  interval $I_{a}$ on $C_{a}$, which is mapped from $I_{b}$ by the inverse
  conformal transformation $\gamma^{(-1)}(z_{b})$.
\end{lemma}

\noindent
Since $f_{b}(\lambda)$ is integrable on $I_{b}$ we have that

\begin{equation}
  \int_{I_{b}}|d\lambda|\,
  |f_{b}(\lambda)|
\end{equation}

\noindent
exists and is finite. If we now consider the integral

\begin{equation}
  \int_{I_{a}}|d\theta|\,
  |f_{a}(\theta)|,
\end{equation}

\noindent
and transform variables from $\theta$ to $\lambda$, recalling that
$f_{a}(\theta)=f_{b}(\lambda)$, we get

\begin{equation}\label{EQAbsInt3}
  \int_{I_{a}}|d\theta|\,
  |f_{a}(\theta)|
  =
  \int_{I_{b}}|d\lambda|\,
  \left|
    \frac{d\theta}{d\lambda}
  \right|
  |f_{b}(\lambda)|.
\end{equation}

\noindent
The absolute value of the derivative shown exists and is finite on
$I_{b}$, given that

\begin{equation}
  \left|
    \frac{d\theta}{d\lambda}
  \right|
  =
  \left|
    \frac{d\gamma^{(-1)}(z_{b})}{dz_{b}}
  \right|,
\end{equation}

\noindent
where $\gamma^{(-1)}(z_{b})$ is analytic on $C_{b}$ and therefore
differentiable there. We also have that $f_{b}(\lambda)$ is integrable on
$I_{b}$. It follows that, since the integrand in the right-hand side of
Equation~(\ref{EQAbsInt3}) is the product of a limited real function with
an integrable real function, and therefore is itself an integrable real
function, the integral in Equation~(\ref{EQAbsInt3}) exists and is finite,
thus implying that $f_{a}(\theta)$ is integrable on the closed interval
$I_{a}$. This establishes Lemma~\ref{Lemma03}.

\vspace{2.6ex}

\noindent
As an immediate consequence of this preliminary result, under the
conditions that we have here, the hypothesis that $f_{b}(\lambda)$ is
locally integrable almost everywhere on $C_{b}$, with the exclusion of the
finite set of points $z_{b,j}$, implies that $f_{a}(\theta)$ is locally
integrable almost everywhere on $C_{a}$, with the exclusion of the
corresponding finite set of points $z_{a,j}$.

\vspace{2.6ex}

\noindent
We must now discuss the issue of the degrees of hardness of the
non-integrable hard singularities of the function $f_{a}(\theta)$ on
$C_{a}$. Since by hypothesis $f_{b}(\lambda)$ has non-integrable hard
singularities at the points $z_{b,j}$, it clearly follows that
$f_{a}(\theta)$ also has hard singularities at the corresponding points
$z_{a,j}$, which may be non-integrable ones. In order to discuss their
degrees of hardness we will use the technique of piecewise integration,
that is, we will consider sectional integrals of $f_{a}(\theta)$ on closed
intervals contained within a neighborhood of the point $z_{a,j}$ where it
has a single isolated hard singularity. Let us show the following
preliminary fact about a real function $f_{a}(\theta)$ defined from a real
function $f_{b}(\lambda)$ which is locally integrable almost everywhere on
$C_{b}$, and which has non-integrable hard singularities with finite
degrees of hardness at the finite set of points $z_{b,j}$.

\begin{lemma}\Colon\label{Lemma04}
  Given a real function $f_{b}(\lambda)$ which has an isolated
  non-integrable hard singularity with finite degree of hardness at a
  point $z_{b,j}$ on $C_{b}$, it follows that the corresponding function
  $f_{a}(\theta)$ defined on the unit circle $C_{a}$ by
  $f_{a}(\theta)=f_{b}(\lambda)$ has an isolated non-integrable hard
  singularity with finite degree of hardness at the corresponding point
  $z_{a,j}$ on $C_{a}$.
\end{lemma}

\noindent
Since the real functions must diverge to infinity at non-integrable hard
singular points, the fact that $f_{a}(\theta)$ has an isolated hard
singularity on $z_{a,j}$ is immediate. Since these singularities are all
isolated from each other, there is on $C_{b}$ a neighborhood of the point
$z_{b,j}$ within which there are no other non-integrable hard
singularities of $f_{b}(\lambda)$. Since the conformal mapping is
continuous, it follows that there is on $C_{a}$ a neighborhood of the
corresponding point $z_{a,j}$ within which there are no other hard
singularities of $f_{a}(\theta)$. Given that the point $z_{a,j}$
corresponds to the angle $\theta_{j}$, let the closed interval
$[\theta_{\ominus,j},\theta_{\oplus,j}]$ contain the point $z_{a,j}$ and
be contained in this neighborhood, so that we have

\begin{equation}
  \theta_{\ominus,j}<\theta_{j}<\theta_{\oplus,j},
\end{equation}

\noindent
where the sole hard singularity of $f_{a}(\theta)$ which is contained
within this interval is the one at the point $\theta_{j}$. Let us now
consider a pair of closed intervals contained within this neighborhood,
one to the left and one to the right of the point $\theta_{j}$, so that we
have

\noindent
\begin{eqnarray}
  I_{\ominus,j}
  & = &
  [\theta_{\ominus,j},\theta_{j}-\varepsilon_{\ominus,j}],
  \nonumber\\
  I_{\oplus,j}
  & = &
  [\theta_{j}+\varepsilon_{\oplus,j},\theta_{\oplus,j}],
\end{eqnarray}

\noindent
where $\varepsilon_{\ominus,j}$ and $\varepsilon_{\oplus,j}$ are two
sufficiently small positive real numbers, so that we also have

\noindent
\begin{eqnarray}
  \theta_{\ominus,j}
  & < &
  \theta_{j}-\varepsilon_{\ominus,j},
  \nonumber\\
  \theta_{j}+\varepsilon_{\oplus,j}
  & < &
  \theta_{\oplus,j}.
\end{eqnarray}

\noindent
Let us now consider sectional primitives of the real function
$f_{a}(\theta)$ on these two intervals. Since the singularity of
$f_{a}(\theta)$ at $\theta_{j}$ may not be integrable, we cannot integrate
across the singularity, but we may integrate within these two lateral
closed intervals, thus defining two sectional primitives of
$f_{a}(\theta)$, one to the left and another one to the right of
$\theta_{j}$,

\noindent
\begin{eqnarray}\label{EQSecInt1}
  f_{a,\ominus}^{-1\prime}(\theta)
  & = &
  \int_{\theta_{0,\ominus,j}}^{\theta}
  d\theta_{\ominus}\,
  f_{a}(\theta_{\ominus}),
  \nonumber\\
  f_{a,\oplus}^{-1\prime}(\theta)
  & = &
  \int_{\theta_{0,\oplus,j}}^{\theta}
  d\theta_{\oplus}\,
  f_{a}(\theta_{\oplus}),
\end{eqnarray}

\noindent
where $f_{a}^{-1\prime}(\theta)$ is the notation for a primitive of
$f_{a}(\theta)$ with respect to $\theta$, where $\theta_{0,\ominus,j}$ and
$\theta_{0,\oplus,j}$ are two arbitrary reference points, one within each
of the two lateral closed intervals, and where we have

\noindent
\begin{equation}
  %
  % Redefining a length: space between rows.
  \renewcommand{\arraystretch}{2.0}
  \begin{array}{rcccl}
    \theta_{\ominus,j}
    & \leq &
    \theta_{\ominus}
    & \leq &
    \theta_{j}-\varepsilon_{\ominus,j},
    \nonumber\\
    \theta_{\ominus,j}
    & \leq &
    \theta_{0,\ominus,j}
    & \leq &
    \theta_{j}-\varepsilon_{\ominus,j},
    \nonumber\\
    \theta_{j}+\varepsilon_{\oplus,j}
    & \leq &
    \theta_{\oplus}
    & \leq &
    \theta_{\oplus,j},
    \nonumber\\
    \theta_{j}+\varepsilon_{\oplus,j}
    & \leq &
    \theta_{0,\oplus,j}
    & \leq &
    \theta_{\oplus,j}.
  \end{array}
\end{equation}

\noindent
If we change variables from $\theta$ to $\lambda$ on the two sectional
integrals in Equation~(\ref{EQSecInt1}), we get

\noindent
\begin{eqnarray}\label{EQSecInt2}
  f_{a,\ominus}^{-1\prime}(\theta)
  & = &
  \int_{\lambda_{0,\ominus,j}}^{\lambda}
  d\lambda_{\ominus}\,
  \frac{d\theta}{d\lambda}(\lambda_{\ominus})\,
  f_{b}(\lambda_{\ominus}),
  \nonumber\\
  f_{a,\oplus}^{-1\prime}(\theta)
  & = &
  \int_{\lambda_{0,\oplus,j}}^{\lambda}
  d\lambda_{\oplus}\,
  \frac{d\theta}{d\lambda}(\lambda_{\oplus})\,
  f_{b}(\lambda_{\oplus}),
\end{eqnarray}

\noindent
where $\lambda_{0,\ominus,j}$ and $\lambda_{0,\oplus,j}$ are the reference
points on $C_{b}$ corresponding respectively to $\theta_{0,\ominus,j}$ and
$\theta_{0,\oplus,j}$. Since the derivative $d\theta/d\lambda$ which
appears in these integrals is finite everywhere on $C_{b}$, and therefore
limited, due to the fact that the inverse conformal transformation
$\gamma^{(-1)}(z_{b})$ is analytic and hence differentiable on the curve
$C_{b}$, it follows that there are two real numbers $R_{m}$ and $R_{M}$
such that

\begin{equation}
  R_{m}
  \leq
  \frac{d\theta}{d\lambda}(\lambda)
  \leq
  R_{M}
\end{equation}

\noindent
everywhere on $C_{b}$. Note that, since the conformal transformation,
besides being continuous and differentiable, is also invertible on
$C_{b}$, the derivative above cannot change sign and thus cannot be zero.
Therefore, the two bounds $R_{m}$ and $R_{M}$ may be chose to have the
same sign. By exchanging the derivative by these extreme values we can
obtain upper and lower bounds for the sectional integrals, and therefore
we get for the sectional primitives in Equation~(\ref{EQSecInt2}),

\noindent
\begin{equation}
  %
  % Redefining a length: space between rows.
  \renewcommand{\arraystretch}{2.0}
  \begin{array}{rcccl}
    R_{m}
    \Int_{\lambda_{0,\ominus,j}}^{\lambda}
    d\lambda_{\ominus}\,
    f_{b}(\lambda_{\ominus})
    & \leq &
    f_{a,\ominus}^{-1\prime}(\theta)
    & \leq &
    R_{M}
    \Int_{\lambda_{0,\ominus,j}}^{\lambda}
    d\lambda_{\ominus}\,
    f_{b}(\lambda_{\ominus}),
    \nonumber\\
    R_{m}
    \Int_{\lambda_{0,\oplus,j}}^{\lambda}
    d\lambda_{\oplus}\,
    f_{b}(\lambda_{\oplus})
    & \leq &
    f_{a,\oplus}^{-1\prime}(\theta)
    & \leq &
    R_{M}
    \Int_{\lambda_{0,\oplus,j}}^{\lambda}
    d\lambda_{\oplus}\,
    f_{b}(\lambda_{\oplus}).
  \end{array}
\end{equation}

\noindent
We now recognize the integrals that appear in these expressions as the
sectional primitives of the function $f_{b}(\lambda)$, so that we get

\noindent
\begin{equation}
  %
  % Redefining a length: space between rows.
  \renewcommand{\arraystretch}{2.0}
  \begin{array}{rcccl}
    R_{m}\,
    f_{b,\ominus}^{-1\prime}(\lambda)
    & \leq &
    f_{a,\ominus}^{-1\prime}(\theta)
    & \leq &
    R_{M}\,
    f_{b,\ominus}^{-1\prime}(\lambda),
    \nonumber\\
    R_{m}\,
    f_{b,\oplus}^{-1\prime}(\lambda)
    & \leq &
    f_{a,\oplus}^{-1\prime}(\theta)
    & \leq &
    R_{M}\,
    f_{b,\oplus}^{-1\prime}(\lambda).
  \end{array}
\end{equation}

\noindent
These expressions are true so long as $\varepsilon_{\ominus,j}$ and
$\varepsilon_{\oplus,j}$, as well as the corresponding quantities
$\delta_{\ominus,j}$ and $\delta_{\oplus,j}$ on $C_{b}$, are not zero, but
since the singularities at $z_{b,j}$ are non-integrable hard ones, we
cannot take the limit in which these quantities tend to zero. Note that,
since $R_{m}$ and $R_{M}$ have the same sign, if the sectional primitives
of the function $f_{b}(\lambda)$ diverge to infinity in this limit, then
so do the sectional primitives of the function $f_{a}(\theta)$. Therefore,
we may conclude that the hard singularity of $f_{a}(\theta)$ is also a
non-integrable one. Since we may take further sectional integrals of these
expressions, without affecting the inequalities, it is immediately
apparent that, after a total of $n$ successive piecewise integrations, we
get for the $n^{\rm th}$ sectional primitives

\noindent
\begin{equation}
  %
  % Redefining a length: space between rows.
  \renewcommand{\arraystretch}{2.0}
  \begin{array}{rcccl}
    R_{m}\,
    f_{b,\ominus}^{-n\prime}(\lambda)
    & \leq &
    f_{a,\ominus}^{-n\prime}(\theta)
    & \leq &
    R_{M}\,
    f_{b,\ominus}^{-n\prime}(\lambda),
    \nonumber\\
    R_{m}\,
    f_{b,\oplus}^{-n\prime}(\lambda)
    & \leq &
    f_{a,\oplus}^{-n\prime}(\theta)
    & \leq &
    R_{M}\,
    f_{b,\oplus}^{-n\prime}(\lambda).
  \end{array}
\end{equation}

\noindent
Since the non-integrable hard singularity of $f_{b}(\lambda)$ at the point
$\lambda_{j}$ which corresponds to $\theta_{j}$ has a finite degree of
hardness, according to the definition of the degrees of hardness, which
was given in~\cite{CAoRFI} and discussed in detail for the case of real
functions in~\cite{CAoRFIV}, there is a value of $n$ such that the limit
in which $\delta_{\ominus,j}\to 0$ and $\delta_{\oplus,j}\to 0$ can be
taken for the sectional primitives $f_{b,\ominus}^{-n\prime}(\lambda)$ and
$f_{b,\oplus}^{-n\prime}(\lambda)$, thus implying that the $n^{\rm th}$
piecewise primitive $f_{b}^{-n\prime}(\lambda)$ of $f_{b}(\lambda)$ is an
integrable real function on the whole interval
$[\lambda_{\ominus,j},\lambda_{\oplus,j}]$ that corresponds to
$[\theta_{\ominus,j},\theta_{\oplus,j}]$, with a borderline hard
singularity, with degree of hardness zero, at the point $\lambda_{j}$. It
follows from the inequalities, therefore, that the corresponding limit in
which $\varepsilon_{\ominus,j}\to 0$ and $\varepsilon_{\oplus,j}\to 0$ can
be taken for the functions $f_{a,\ominus}^{-n\prime}(\theta)$ and
$f_{a,\oplus}^{-n\prime}(\theta)$, thus implying that the $n^{\rm th}$
piecewise primitive $f_{a}^{-n\prime}(\theta)$ of $f_{a}(\theta)$ is also
an integrable real function on the whole interval
$[\theta_{\ominus,j},\theta_{\oplus,j}]$, with a borderline hard
singularity, with degree of hardness zero, at $\theta_{j}$. Therefore, the
non-integrable hard singularity of $f_{a}(\theta)$ at $\theta_{j}$ has a
finite degree of hardness, to wit the same degree of hardness $n$ of the
corresponding non-integrable hard singularity of $f_{b}(\lambda)$ at
$\lambda_{j}$. This establishes Lemma~\ref{Lemma04}.

\vspace{2.6ex}

\noindent
We have therefore established that, so long as $f_{b}(\lambda)$ is locally
integrable almost everywhere on $C_{b}$, and so long as its non-integrable
hard singularities have finite degrees of hardness, these same two facts
are true for $f_{a}(\theta)$ on $C_{a}$. Since we thus see that the
necessary properties of the real functions are preserved by the conformal
transformation, we are therefore in a position to use the result of
Theorem~\ref{Theo02} in Section~\ref{Sec03} in order to extend the
existence theorem of the Dirichlet problem to non-integrable real
functions which are, however, integrable almost everywhere on $C_{b}$,
still for the case of differentiable curves.

\vspace{2.6ex}

\noindent
In this section, using once again the results from the previous sections,
we will establish the following theorem.

\begin{theorem}\Colon\label{Theo05}
  Given a differentiable simple closed curve $C$ of finite total length on
  the complex plane $z=x+\ii y$, given the invertible conformal
  transformation $\gamma(z)$ whose derivative has no zeros on the closed
  unit disk, that maps it from the unit circle, and given a real function
  $f(\lambda)$ on that curve, that satisfies the list of conditions
  described below, there is a solution $u(x,y)$ of the Dirichlet problem
  of the Laplace equation within the interior of that curve, that assumes
  the given values $f(\lambda)$ almost everywhere at the curve.
\end{theorem}

\begin{proof}\Colon
\end{proof}

\noindent
Similarly to what was done in the two previous cases, in
Sections~\ref{Sec05} and~\ref{Sec06}, the proof consists of using the
conformal transformation between the closed unit disk and the union of the
curve $C$ with its interior, which according to the analysis in
Section~\ref{Sec04} always exists, to map the given boundary condition on
$C$ onto a corresponding boundary condition on the unit circle, then using
the proof of existence established before by Theorem~\ref{Theo02} in
Section~\ref{Sec03} for the closed unit disk to establish the existence of
the solution of the corresponding Dirichlet problem on that disk, and
finally using once more the conformal transformation to map the resulting
solution back to $C$ and its interior, thus obtaining the solution of the
original Dirichlet problem. The list of conditions on the real functions
is now the following.

\begin{enumerate}

\item The real function $f(\lambda)$ is locally integrable almost
  everywhere on $C$, including the cases in which this function is
  globally integrable there.

\item The number of hard singularities on the unit circle of the
  corresponding inner analytic function $w(z)$ on the unit disk is finite.

\item The hard singularities of the corresponding inner analytic function
  $w(z)$ have finite degrees of hardness.

\end{enumerate}

\noindent
The rest of the proof is identical to that of the two previous cases, in
Sections~\ref{Sec05} and~\ref{Sec06}. Therefore, once again we may
conclude that, due to the existence theorem of the Dirichlet problem on
the unit disk of the plane $z_{a}$, which in this case was established by
Theorem~\ref{Theo02} in Section~\ref{Sec03}, we know that there is an
inner analytic function $w_{a}(z_{a})$ such that its real part
$u_{a}(\rho,\theta)$ is harmonic within the open unit disk and satisfies
$u_{a}(1,\theta)=f_{a}(\theta)$ almost everywhere at the boundary
$C_{a}$. Just as before, we get on the $z_{b}$ plane the complex function
$w_{b}(z_{b})$ which is analytic in the interior of the curve
$C_{b}$. Therefore, the real part $u_{b}(x,y)$ of $w_{b}(z_{b})$ is
harmonic and thus satisfies

\begin{equation}
  \nabla^{2}u_{b}(x,y)
  =
  0
\end{equation}

\noindent
in the interior of $C_{b}$, while we also have that

\begin{equation}
  u_{b}(x,y)=f_{b}(\lambda),
\end{equation}

\noindent
almost everywhere on $C_{b}$. This establishes the existence, by
construction, of the solution of the Dirichlet problem on the $z_{b}$
plane, under our current hypotheses. This completes the proof of
Theorem~\ref{Theo05}.

\vspace{2.6ex}

\noindent
In this way we have generalized the proof of existence of the Dirichlet
problem from the unit circle to all differentiable simple closed curves
with finite total lengths on the plane, for boundary conditions given by
non-integrable real functions which are locally integrable almost
everywhere and have at most a finite set of hard singular points.

\section{Non-Integrable Functions on Non-Differentiable
  Curves}\label{Sec08}

In this section we will show how one can generalize the proof of existence
of the Dirichlet problem on the unit disk, given by Theorem~\ref{Theo02}
in Section~\ref{Sec03}, to the case in which we have, as the boundary
condition, non-integrable real functions $f_{b}(\lambda)$ defined on
boundaries given by non-differentiable curves $C_{b}$ on the plane. We
will be able to do this if the non-integrable real functions, despite
being non-integrable over the whole curves $C_{b}$, are locally integrable
almost everywhere on those curves, and if, in addition to this, the
non-integrable hard singularities involved have finite degrees of
hardness. The definition of the concept of local integrability almost
everywhere is that given by Definition~\ref{Def04}, in
Section~\ref{Sec07}. The additional difficulty that appears in this case
is the same which was discussed in Section~\ref{Sec06}, due to the fact
that the curve $C_{b}$ is not differentiable at the finite set of points
$z_{b,i}$, for $i\in\{1,\ldots,N\}$. Just as in that case, this will
require that we impose one additional limitation on the real functions
giving the boundary conditions, namely that any integrable hard
singularities where they diverge to infinity do not coincide with any of
the points $z_{b,i}$.

The proof will follow the general lines of the one given for integrable
real functions in Section~\ref{Sec06}, with the difference that, just as
we did in Section~\ref{Sec07}, instead of showing that the corresponding
functions $f_{a}(\theta)$ on the unit circle $C_{a}$ are integrable there,
we will show that they are locally integrable almost everywhere there. In
addition to this, instead of using the existence theorem for integrable
real function on the unit circle, which was given by Theorem~\ref{Theo01}
in Section~\ref{Sec02}, we will use the corresponding result for
non-integrable real functions which are locally integrable almost
everywhere on the unit circle, which was given by Theorem~\ref{Theo02} in
Section~\ref{Sec03}. Since that result depends also on the non-integrable
hard singularities of the real functions having finite degrees of
hardness, we will also show that the hypothesis that the functions
$f_{b}(\lambda)$ have this property implies that the corresponding
functions $f_{a}(\theta)$ have the same property as well. In order to do
this we will use once again the technique of piecewise integration which
was introduced in~\cite{CAoRFIV}.

The preliminary result given by Lemma~\ref{Lemma03} in Section~\ref{Sec07}
is still valid here. As a consequence of this we may conclude at once
that, under the conditions that we have here, the hypothesis that
$f_{b}(\lambda)$ is locally integrable almost everywhere on $C_{b}$, with
the exclusion of the finite set of points $z_{b,j}$, for
$j\in\{1,\ldots,M\}$, implies that $f_{a}(\theta)$ is locally integrable
almost everywhere on the unit circle $C_{a}$, with the exclusion of the
corresponding finite set of points $z_{a,j}=\gamma^{(-1)}(z_{b,j})$.

We must now discuss the issue of the degrees of hardness of the
non-integrable hard singularities of the function $f_{a}(\theta)$ on
$C_{a}$. Since by hypothesis $f_{b}(\lambda)$ has non-integrable hard
singularities at the points $z_{b,j}$, it clearly follows that
$f_{a}(\theta)$ also has hard singularities at the corresponding points
$z_{a,j}$, which may be non-integrable ones. In order to discuss their
degrees of hardness we will use the technique of piecewise integration,
that is, we will consider sectional integrals of $f_{a}(\theta)$ on closed
intervals contained within a neighborhood of the point $z_{a,j}$ where it
has a single isolated hard singularity. Let us show the following
preliminary fact about a real function $f_{a}(\theta)$ defined from a real
function $f_{b}(\lambda)$ which is locally integrable almost everywhere on
$C_{b}$, and which has non-integrable hard singularities at the finite set
of points $z_{b,j}$.

\begin{lemma}\Colon\label{Lemma05}
  Given a real function $f_{b}(\lambda)$ which has an isolated
  non-integrable hard singularity with finite degree of hardness at a
  point $z_{b,j}$ on $C_{b}$, and which is such that the hard
  singularities where it diverges to infinity do not coincide with any of
  the points $z_{b,i}$ where $C_{b}$ is non-differentiable, it follows
  that the corresponding function $f_{a}(\theta)$ defined on the unit
  circle $C_{a}$ by $f_{a}(\theta)=f_{b}(\lambda)$ has an isolated
  non-integrable hard singularity with finite degree of hardness at the
  corresponding point $z_{a,j}$ on $C_{a}$.
\end{lemma}

\noindent
Since the real functions must diverge to infinity at non-integrable hard
singular points, the fact that $f_{a}(\theta)$ has an isolated hard
singularity on $z_{a,j}$ is immediate. Since these singularities are all
isolated from each other, and since they do not coincide with the points
$z_{b,i}$ where $C_{b}$ is non-differentiable, there is on $C_{b}$ a
neighborhood of the point $z_{b,j}$ within which $C_{b}$ is differentiable
and there are no other non-integrable hard singularities of
$f_{b}(\lambda)$. Since the conformal mapping is continuous, it follows
that there is on $C_{a}$ a neighborhood of the corresponding point
$z_{a,j}$ within which there are no zeros of the derivative of
$\gamma(z_{a})$ and no other hard singularities of $f_{a}(\theta)$. The
construction of the two sectional primitives of $f_{a}(\theta)$ by means
of piecewise integration is the same as the one which was executed before
for Lemma~\ref{Lemma04} in Section~\ref{Sec07}, resulting in
Equation~(\ref{EQSecInt1}). After we change variables from $\theta$ to
$\lambda$ on the two sectional integrals in that equation we get

\noindent
\begin{eqnarray}\label{EQSecInt3}
  f_{a,\ominus}^{-1\prime}(\theta)
  & = &
  \int_{\lambda_{0,\ominus,j}}^{\lambda}
  d\lambda_{\ominus}\,
  \frac{d\theta}{d\lambda}(\lambda_{\ominus})\,
  f_{b}(\lambda_{\ominus}),
  \nonumber\\
  f_{a,\oplus}^{-1\prime}(\theta)
  & = &
  \int_{\lambda_{0,\oplus,j}}^{\lambda}
  d\lambda_{\oplus}\,
  \frac{d\theta}{d\lambda}(\lambda_{\oplus})\,
  f_{b}(\lambda_{\oplus}),
\end{eqnarray}

\noindent
where all the symbols involved are the same as before. Since the
derivative $d\theta/d\lambda$ which appears in these integrals is finite,
and therefore limited, everywhere within the two lateral intervals
involved, due to the fact that the inverse conformal transformation
$\gamma^{(-1)}(z_{b})$ is analytic and therefore differentiable within
these intervals, it follows that there are two pairs of real numbers
$R_{\ominus,m}$ and $R_{\ominus,M}$, as well as $R_{\oplus,m}$ and
$R_{\oplus,M}$, such that

\noindent
\begin{equation}
  %
  % Redefining a length: space between rows.
  \renewcommand{\arraystretch}{2.0}
  \begin{array}{rcccl}
    R_{\ominus,m}
    & \leq &
    \FFrac{d\theta}{d\lambda}(\lambda_{\ominus})
    & \leq &
    R_{\ominus,M},
    \nonumber\\
    R_{\oplus,m}
    & \leq &
    \FFrac{d\theta}{d\lambda}(\lambda_{\oplus})
    & \leq &
    R_{\oplus,M},
  \end{array}
\end{equation}

\noindent
everywhere within each interval. Note that, since the conformal
transformation, besides being continuous and differentiable, is also
invertible within each interval, the derivatives above cannot change sign
and thus cannot be zero. Therefore, the pair of bounds $R_{\ominus,m}$ and
$R_{\ominus,M}$ may be chosen to have the same sign, and so may the pair
of bounds $R_{\oplus,m}$ and $R_{\oplus,M}$. By exchanging the derivative
by these extreme values we can obtain upper and lower bounds for the
sectional integrals, and therefore we get for the sectional primitives in
Equation~(\ref{EQSecInt3}),

\noindent
\begin{equation}
  %
  % Redefining a length: space between rows.
  \renewcommand{\arraystretch}{2.0}
  \begin{array}{rcccl}
    R_{\ominus,m}
    \Int_{\lambda_{0,\ominus,j}}^{\lambda}
    d\lambda_{\ominus}\,
    f_{b}(\lambda_{\ominus})
    & \leq &
    f_{a,\ominus}^{-1\prime}(\theta)
    & \leq &
    R_{\ominus,M}
    \Int_{\lambda_{0,\ominus,j}}^{\lambda}
    d\lambda_{\ominus}\,
    f_{b}(\lambda_{\ominus}),
    \nonumber\\
    R_{\oplus,m}
    \Int_{\lambda_{0,\oplus,j}}^{\lambda}
    d\lambda_{\oplus}\,
    f_{b}(\lambda_{\oplus})
    & \leq &
    f_{a,\oplus}^{-1\prime}(\theta)
    & \leq &
    R_{\oplus,M}
    \Int_{\lambda_{0,\oplus,j}}^{\lambda}
    d\lambda_{\oplus}\,
    f_{b}(\lambda_{\oplus}).
  \end{array}
\end{equation}

\noindent
We now recognize the integrals that appear in these expressions as the
sectional primitives of the function $f_{b}(\lambda)$, so that we get

\noindent
\begin{equation}
  %
  % Redefining a length: space between rows.
  \renewcommand{\arraystretch}{2.0}
  \begin{array}{rcccl}
    R_{\ominus,m}\,
    f_{b,\ominus}^{-1\prime}(\lambda)
    & \leq &
    f_{a,\ominus}^{-1\prime}(\theta)
    & \leq &
    R_{\ominus,M}\,
    f_{b,\ominus}^{-1\prime}(\lambda),
    \nonumber\\
    R_{\oplus,m}\,
    f_{b,\oplus}^{-1\prime}(\lambda)
    & \leq &
    f_{a,\oplus}^{-1\prime}(\theta)
    & \leq &
    R_{\oplus,M}\,
    f_{b,\oplus}^{-1\prime}(\lambda).
  \end{array}
\end{equation}

\noindent
These expressions are true so long as $\varepsilon_{\ominus,j}$ and
$\varepsilon_{\oplus,j}$, as well as the corresponding quantities
$\delta_{\ominus,j}$ and $\delta_{\oplus,j}$ on $C_{b}$, are not zero, but
since the singularities at $z_{b,j}$ are non-integrable hard ones, we
cannot take the limit in which these quantities tend to zero. Note that,
since $R_{\ominus,m}$ and $R_{\ominus,M}$ have the same sign, and also
$R_{\oplus,m}$ and $R_{\oplus,M}$ have the same sign, if the sectional
primitives of the function $f_{b}(\lambda)$ diverge to infinity in this
limit, then so do the sectional primitives of the function
$f_{a}(\theta)$. Therefore, we may conclude that the hard singularity of
$f_{a}(\theta)$ is also a non-integrable one. Since we may take further
sectional integrals of these expressions, without affecting the
inequalities, it is immediately apparent that, after a total of $n$
successive piecewise integrations, we get for the $n^{\rm th}$ sectional
primitives

\noindent
\begin{equation}
  %
  % Redefining a length: space between rows.
  \renewcommand{\arraystretch}{2.0}
  \begin{array}{rcccl}
    R_{\ominus,m}\,
    f_{b,\ominus}^{-n\prime}(\lambda)
    & \leq &
    f_{a,\ominus}^{-n\prime}(\theta)
    & \leq &
    R_{\ominus,M}\,
    f_{b,\ominus}^{-n\prime}(\lambda),
    \nonumber\\
    R_{\oplus,m}\,
    f_{b,\oplus}^{-n\prime}(\lambda)
    & \leq &
    f_{a,\oplus}^{-n\prime}(\theta)
    & \leq &
    R_{\oplus,M}\,
    f_{b,\oplus}^{-n\prime}(\lambda).
  \end{array}
\end{equation}

\noindent
Since the non-integrable hard singularity of $f_{b}(\lambda)$ at the point
$\lambda_{j}$ which corresponds to $\theta_{j}$ has a finite degree of
hardness, according to the definition of the degrees of hardness, which
was given in~\cite{CAoRFI} and discussed in detail for the case of real
functions in~\cite{CAoRFIV}, there is a value of $n$ such that the limit
in which $\delta_{\ominus,j}\to 0$ and $\delta_{\oplus,j}\to 0$ can be
taken for the sectional primitives $f_{b,\ominus}^{-n\prime}(\lambda)$ and
$f_{b,\oplus}^{-n\prime}(\lambda)$, thus implying that the $n^{\rm th}$
piecewise primitive $f_{b}^{-n\prime}(\lambda)$ of $f_{b}(\lambda)$ is an
integrable real function on the whole interval
$[\lambda_{\ominus,j},\lambda_{\oplus,j}]$, with a borderline hard
singularity, with degree of hardness zero, at the point $\lambda_{j}$. It
follows from the inequalities, therefore, that the corresponding limit in
which $\varepsilon_{\ominus,j}\to 0$ and $\varepsilon_{\oplus,j}\to 0$ can
be taken for the functions $f_{a,\ominus}^{-n\prime}(\theta)$ and
$f_{a,\oplus}^{-n\prime}(\theta)$, thus implying that the $n^{\rm th}$
piecewise primitive $f_{a}^{-n\prime}(\theta)$ of $f_{a}(\theta)$ is also
an integrable real function on the whole interval
$[\theta_{\ominus,j},\theta_{\oplus,j}]$, with a borderline hard
singularity, with degree of hardness zero, at $\theta_{j}$. Therefore, the
non-integrable hard singularity of $f_{a}(\theta)$ at $\theta_{j}$ has a
finite degree of hardness, to wit the same degree of hardness $n$ of the
corresponding non-integrable hard singularity of $f_{b}(\lambda)$ at
$\lambda_{j}$. This establishes Lemma~\ref{Lemma05}.

\vspace{2.6ex}

\noindent
We have therefore established that, under the hypothesis that the hard
singularities where $f_{b}(\lambda)$ diverges to infinity do not coincide
with any of the points $z_{b,i}$ where $C_{b}$ is non-differentiable, so
long as $f_{b}(\lambda)$ is locally integrable almost everywhere on
$C_{b}$, and so long as its non-integrable hard singularities have finite
degrees of hardness, these same two facts are true for $f_{a}(\theta)$ on
$C_{a}$. Since we thus see that the necessary properties of the real
functions are preserved by the conformal transformation, we are therefore
in a position to use the result of Theorem~\ref{Theo02} in
Section~\ref{Sec03} in order to extend the existence theorem of the
Dirichlet problem to non-integrable real functions which are, however,
integrable almost everywhere on $C_{b}$, this time for the case of
non-differentiable curves.

\vspace{2.6ex}

\noindent
In this section, using one more time the results from the previous
sections, we will establish the following theorem.

\begin{theorem}\Colon\label{Theo06}
  Given a simple closed curve $C$ of finite total length on the complex
  plane $z=x+\ii y$, which is non-differentiable at a given finite set of
  points $z_{i}$, for $i\in\{1,\ldots,N\}$, given the conformal
  transformation $\gamma(z)$ that maps it from the unit circle, whose
  derivative has zeros on the unit circle at the corresponding set of
  points, and given a real function $f(\lambda)$ on that curve, that
  satisfies the list of conditions described below, there is a solution
  $u(x,y)$ of the Dirichlet problem of the Laplace equation within the
  interior of that curve, that assumes the given values $f(\lambda)$
  almost everywhere at the curve.
\end{theorem}

\begin{proof}\Colon
\end{proof}

\noindent
Similarly to what was done in the three previous sections, the proof
consists of using the conformal transformation between the closed unit
disk and the union of the curve $C$ with its interior, which according to
the analysis in Section~\ref{Sec04} always exists, to map the given
boundary condition on $C$ onto a corresponding boundary condition on the
unit circle, then using the proof of existence established before by
Theorem~\ref{Theo02} in Section~\ref{Sec03} for the closed unit disk to
establish the existence of the solution of the corresponding Dirichlet
problem on that disk, and finally using once more the conformal
transformation to map the resulting solution back to $C$ and its interior,
thus obtaining the solution of the original Dirichlet problem. The list of
conditions on the real functions is now the following.

\begin{enumerate}

\item The real function $f(\lambda)$ is locally integrable almost
  everywhere on $C$, including the cases in which this function is
  globally integrable there.

\item The number of hard singularities on the unit circle of the
  corresponding inner analytic function $w(z)$ on the unit disk is finite.

\item The hard singularities of the corresponding inner analytic function
  $w(z)$ have finite degrees of hardness.

\item The hard singularities of $f(\lambda)$ where it diverges to infinity
  are not located at any of the points where the curve $C$ is
  non-differentiable.

\end{enumerate}

\noindent
The rest of the proof is identical to that of the three previous cases.
Therefore, once again we may conclude that, due to the existence theorem
of the Dirichlet problem on the unit disk of the plane $z_{a}$, which in
this case was established in by Theorem~\ref{Theo02} Section~\ref{Sec03},
we know that there is an inner analytic function $w_{a}(z_{a})$ such that
its real part $u_{a}(\rho,\theta)$ is harmonic within the open unit disk
and satisfies $u_{a}(1,\theta)=f_{a}(\theta)$ almost everywhere at the
boundary $C_{a}$. Just as before, we get on the $z_{b}$ plane the complex
function $w_{b}(z_{b})$ which is analytic in the interior of the curve
$C_{b}$. Therefore, the real part $u_{b}(x,y)$ of $w_{b}(z_{b})$ is
harmonic and thus satisfies

\begin{equation}
  \nabla^{2}u_{b}(x,y)
  =
  0
\end{equation}

\noindent
in the interior of $C_{b}$, while we also have that

\begin{equation}
  u_{b}(x,y)=f_{b}(\lambda),
\end{equation}

\noindent
almost everywhere on $C_{b}$. This establishes the existence, by
construction, of the solution of the Dirichlet problem on the $z_{b}$
plane, under our current hypotheses. This completes the proof of
Theorem~\ref{Theo06}.

\vspace{2.6ex}

\noindent
In this way we have generalized the proof of existence of the Dirichlet
problem from the unit circle to all simple closed curves with finite total
lengths on the plane, that can be either differentiable or
non-differentiable on a finite set of points, but now for boundary
conditions given by non-integrable real functions which are locally
integrable almost everywhere and have at most a finite set of hard
singular points.

\section{Conclusions and Outlook}\label{Sec09}

A very general proof of the existence of the solution of the Dirichlet
boundary value problem of the Laplace equation on the plane was presented.
The proof is valid not only for a very large class of real functions at
the boundary, but also for a large class of boundary curves, with and
without points of non-differentiability. The proof was presented in
incremental steps, each generalizing the previous ones. The proofs for the
unit circle are based on the complex-analytic structure within the unit
disk presented and developed in previous
papers~\cite{CAoRFI,CAoRFII,CAoRFIII,CAoRFIV}. The generalization for
curves other than the unit circle uses the conformal mapping results
associated to the famous Riemann mapping theorem. The most general
statement of the theorem established here reads as follows.

\vspace{3ex}

\noindent
\parbox{\textwidth}{\bf Given a real function $f(\lambda)$ that defines
  the boundary condition on a plane curve $C$ parametrized by the real
  arc-length variable $\lambda$, so long as the real function is locally
  integrable almost everywhere on $C$, and is such that the corresponding
  inner analytic function has at most a finite number of hard
  singularities with finite degrees of hardness, so long as $C$ is a
  simple closed curve with finite total length and at most a finite number
  of points of non-differentiability, and so long as the hard singular
  points of $f(\lambda)$ where it diverges to infinity do not coincide
  with any of the points where the curve is not differentiable, there
  exists a real function $u(x,y)$ that satisfies $\nabla^{2}u(x,y)=0$ in
  the interior of $C$ and that satisfies $u(x,y)=f(\lambda)$ almost
  everywhere on $C$.}

\vspace{3ex}

\noindent
The proof is constructive, and consists of constructing from $f(\lambda)$
an analytic function in the interior of $C$, of which $u(x,y)$ is the real
part. The theorem is quite general, including large classes of both
boundary conditions and boundary curves.

Further extensions of the theorem may be possible. For example, the proofs
established in Sections~\ref{Sec02} and~\ref{Sec03} can be rather
trivially extended to include as well the whole space of singular Schwartz
distributions discussed in~\cite{CAoRFII}, that is, they can be extended
to generalized real functions. This allows one to discuss some rather
unusual Dirichlet problems in which the boundary condition is given by a
singular real object such as the Dirac delta ``function'' or its
derivatives. As mentioned in~\cite{CAoRFI}, a possible further extension
would be to real functions with a countable infinity of hard singular
points which have, however, a finite number of accumulation points. The
requirement that the hard singular points of $f(\lambda)$ where it
diverges to infinity do not coincide with the points where the curve $C$
is non-differentiable seems to be a technical quirk, and probably can be
eliminated. It is important to note that the proof is intrinsically
limited to two-dimensional problems on the plane.

It is interesting to observe that the uniqueness of the solution can also
be discussed in this context, in terms of the fact that $f(\theta)\equiv
0$ corresponds to the Fourier coefficients $\alpha_{0}=0$, $\alpha_{k}=0$
and $\beta_{k}=0$, for all $k$, and therefore to the complex Taylor
coefficients $c_{0}=0$ and $c_{k}=0$, for all $k$, and therefore to the
identically zero inner analytic function $w(z)\equiv 0$. Given an
integrable real function $f(\theta)$ and two corresponding solutions
$w_{1}(z)$ and $w_{2}(z)$ of the Dirichlet problem, we simply consider
$w(z)=w_{2}(z)-w_{1}(z)$, which is therefore a solution of the Dirichlet
problem with $f(\theta)\equiv 0$, and thus by construction is $w(z)\equiv
0$. It follows that $w_{2}(z)\equiv w_{1}(z)$, so that the solution is
unique, in the sense that $u_{2}(\rho,\theta)=u_{1}(\rho,\theta)$ almost
everywhere on the unit disk. We can say, in fact, that these two functions
are equal at all points on the unit circle where they are well defined.

With some more work towards its generalization, the result presented here
points, perhaps, to an even more general result, according to which the
solution of the Dirichlet problem of the Laplace equation in two
dimensions, in essence, {\em always} exists, in the sense that it exists
under all conceivable circumstances in which it makes any sense at all to
pose the corresponding boundary value problem. Already, even with the
result as it is now, this is almost the case in what concerns the
applications to Physics.

\section*{Acknowledgments}

The author would like to thank his friend and colleague Prof. Carlos
Eugênio Imbassay Carneiro, to whom he is deeply indebted for all his
interest and help, as well as his careful reading of the manuscript and
helpful criticism regarding this work.

\bibliography{allrefs_en}\bibliographystyle{ieeetr}

\end{document}

%% file: packages.tex
\usepackage[latin1]{inputenc}
\usepackage{epsfig}
\usepackage{amsfonts}
\usepackage{cite}
\usepackage{color}
\definecolor{lgrey}{gray}{0.99}

%% file: math-defs.tex
\newcommand{\ii}{\mbox{\boldmath$\imath$}}

\newcommand{\FFrac}[2]{{\displaystyle\frac{\displaystyle #1}{\displaystyle #2}}}

\newcommand{\Int}{\displaystyle\int}

\newcommand{\ldot}{\mbox{\Large$\cdot$}\!}

\newtheorem{definition}{Definition}

\newtheorem{lemma}{Lemma}
\newtheorem{theorem}{Theorem}

\newtheorem{proof}{Proof}[theorem]
\newcommand{\Colon}{{\hspace{-0.4em}\bf:}}